\documentclass[smallextended]{svjour3}       
\smartqed  
\usepackage{graphicx}
\usepackage{subfigure}
\usepackage{amsmath}
\usepackage{tabu}
\usepackage[numbers]{natbib}
\begin{document}

\title{Insights on aliasing driven instabilities for advection equations with application to Gauss-Lobatto discontinuous Galerkin methods
}

\titlerunning{Aliasing driven instabilities on DGSEM-GL methods}        

\author{Juan Manzanero         \and
        Gonzalo Rubio \and 
        Esteban Ferrer \and 
        Eusebio Valero \and \\
        David A. Kopriva 
}


\institute{J. Manzanero (\email{juan.manzanero@upm.es}) \and G. Rubio \and E. Ferrer \and E. Valero \at
              ETSIAE-UPM - Universidad Polit\'ecnica de Madrid. School of Aeronautics, Plaza Cardenal Cisneros 3, E-28040 Madrid, Spain. \\
			\and 
			David A. Kopriva \at Department of Mathematics, The Florida State University, Tallahassee, FL 32306, USA.
}

\date{Received: date / Accepted: date}

\maketitle

\begin{abstract}
We analyse instabilities due to aliasing errors when solving one dimensional non-constant advection speed equations and discuss means to alleviate these types of errors when using high order discontinuous Galerkin (DG) schemes.
First, we compare analytical bounds for the continuous and discrete version of the PDEs.
Whilst traditional $L^2$ norm energy bounds applied to the discrete PDE do not always predict the physical behaviour of the continuous version of the equation, more strict elliptic norm bounds correctly bound the behaviour of the continuous PDE.
Having derived consistent bounds, we analyse the effectiveness of two stabilising techniques: over-integration and split form variations (conservative, non-conservative and skew-symmetric). Whilst the former is shown to not alleviate aliasing in general, the latter ensures an aliasing-free solution if the splitting form of the discrete PDE is consistent with the continuous equation. The success of split form de-aliasing is restricted
to DG schemes with the summation-by-parts simultaneous-approximation-term (SBP-SAT) properties (e.g. DG with Gauss-Lobatto points). Numerical experiments are included to illustrate the theoretical findings.
\end{abstract}

\section{Introduction} 



High order methods are the preferred discretisation technique when high accuracy is required \cite{WorkshopDG}. The low numerical dissipation inherent to high-order methods, however, may not be sufficient to mask instabilities such as aliasing errors arising from the non-linear flux discretisation \cite{2006:KirbyAliasing}.
In Spectral/hp element methods, the stabilisation of the scheme has been traditionally achieved by means of over-integration techniques (see \cite{2015:Spiegel}), by using skew symmetric forms \cite{zang1991rotation}, or by adding artificial viscosity \cite{2014:Guermond}. The last can be implemented by upwinding inter-element fluxes (in discontinuous schemes, e.g. DG, see \cite{2009:Toro}), or by including an artificial dissipation term that maintains high-order accuracy (spectral vanishing viscosity methods \cite{2000:Karamanos}).

A different approach, followed by some in the finite differences community, is to develop schemes that are provably stable thanks to the summation-by-parts (SBP) and simultaneous-approximation-term (SAT) properties. The summation-by-parts property can be used to prove stability in each element (see \cite{2013:Fisher}), whilst the simultaneous-approximation-term is fundamental to prove stability at both interior and physical boundaries (see \cite{2013:Carpenter}). A review on SBP-SAT finite differences schemes, and their stability properties can be found in \cite{2014:Nordstrom}. 

Recent research on Spectral/hp methods that satisfy the summation-by-parts property has enabled the adaption of stabilisation techniques and theorems developed in the finite differences community, while retaining favourable high-order spectral element properties. In \cite{2013:Gassner}, for instance, Gassner showed that the discontinuous Galerkin spectral element method with Gauss-Lobatto points satisfies all formal definitions of an SBP-SAT scheme, allowing him to obtain methods that are provably stable. In \cite{2013:Gassner} an energy stable split formulation scheme was developed for the Burgers equation. In \cite{Kopriva2}, schemes that are strongly stable were constructed for the non-constant advection equation, considering multidimensional systems with curved elements. Later, in \cite{2015:kopriva} the methodology was extended to moving geometries, under the ALE (Arbitrary-Lagrangian-Eulerian) formulation. The last step forward made by Gassner et. al. has been the derivation of schemes that are provably stable for the Euler fluid dynamics equations \cite{2016:gassner}, the Magneto-Hydrodynamics \cite{2016:Winters} and the Navier-Stokes equations \cite{gassner2017br1}. 

Here, we extend the work of Kopriva and Gassner \cite{Kopriva2}, where schemes that are strongly stable (consistent with the energy bounds of the continuous PDE) were developed for the non-constant advection equation. The traditional energy bounds used there do not always predict the behaviour of the analytical solution so we use more strict bounds (based on elliptic norms \cite{fornberg1975fourier,gottlieb2001spectral}) to bound correctly the behaviour of the continuous PDE. Numerical schemes consistent with these new bounds produce numerical solutions with the same behaviour as the analytical solution. 
We first derive these bounds for the continuous PDE to then extend them, using the SBP-SAT properties of the DGSEM with Gauss-Lobatto nodes, to the discretisation. Then we analyse the effectiveness of two stabilising techniques: over-integration and split form variations (conservative, non-conservative and skew-symmetric).

The paper is organised as follows: we first present a summary of the results in Section \ref{sec:Summary}. We then introduce the PDE to be studied and its traditional energy bounds in Section \ref{sec:Preliminaries}. In Section \ref{sec:AlternativeContinuousBounds}, we derive alternative energy bounds (based on elliptic norms) for the continuous PDE. Then, in Section \ref{sec:DGSEMversion} we derive a discrete version for the DGSEM with Gauss-Lobatto nodes, with a discussion on stabilising techniques. Next, in Section \ref{sec:numericalExperiments}, we present numerical experiments to show the validity of the stability bounds. Moreover, the Gauss nodes version of the DGSEM method is recovered to show that the capability to perform stabilisation by means of split formulations does not hold with that set of nodes. Lastly, we give some final remarks in Section \ref{sec:conclusions}.

\section{Summary of Results: Continuous and Discrete Energy Bounds}\label{sec:Summary}
Consider the first order one-dimensional initial-boundary-value problem:

\begin{equation}
\begin{split}
u_t + f_x &= F,\quad 0 < x < L, \; t > 0,\\
u(x,0) &= u_0(x),	\\
u(0,t) &= g_L(t), \\
u(L,t) &= g_R(t),
\end{split}
\label{eq:general_advection_problem}
\end{equation}
where $u(x,t)$ is the solution, $f=f(u,x) = a(x)u$ is the flux, and $F=F(u,x)$ is a source term that might depend on the solution itself, and explicitly on the space coordinate. For advection equations, we show, by an alternative bounds technique, how aliasing errors lead to numerical instabilities in the discrete solution. Our findings are summarised in Table \ref{tab:results_summary}. Notice that, $a(x)$ is the advection speed, whilst $I^N(a)$ refers to  its interpolation to the Gauss-Lobatto points. Any exponential growth experienced by the solution energy is entirely driven by aliasing errors, which are measured with the coefficient $\gamma$. This coefficient depends on the advection speed function, $a(x)$ and the polynomial order, $N$, of the approximation. The parameter $\alpha$ controls the discretisation splitting, and must be chosen properly to cancel the aliasing errors to obtain a stable scheme. For the conservative advection equation, the value $\alpha=1$ (i.e. conservative DG) suffices to avoid the energy growth, whilst the value $\alpha=0$ (i.e. non-conservative DG) suffices to remove any aliasing-driven energy growth when approximating the non-conservative equation.

\begin{table}[h!]
	\centering
	\caption[caption]{Stability of continuous advection PDEs, and their discontinuous Galerkin discrete version. }
	\label{tab:results_summary}	
	
	\bgroup
	\def\arraystretch{1}
	\begin{tabu} to \textwidth{ll}
		\hline \\
		PDE: & Conservative advection equation \\\\ \hline \\
		Equation & $u_t + (au)_x = 0,~~~a(x)>0 ~;~~~ f=a(x)u$ \\\\
		Split version & $u_t + \alpha (au)_x + (1-\alpha)(a_x u + au_x) = 0$\\\\
		& \hspace*{0.5cm}$\bullet$ $\alpha=0$: Non-conservative DG \\\\
		& \hspace*{0.5cm}$\bullet$ $\alpha=1/2$: Skew-symmetric DG \\\\
		& \hspace*{0.5cm}$\bullet$ $\alpha=1$: Conservative DG \\\\		  								
		Continuous bound & $\displaystyle{\min_{x\in[0,L]}\{a(x)\}||u(T)||^2 \leq \max_{x\in[0,L]}\{a(x)\}||u_0||^2 + \int_{0}^T a(0)^2 g_L^2 dt}$\\\\
		Discrete bound  &$\displaystyle{\min_{x\in[0,L]}\{I^N(a)\}||U(T)||^2 \leq 3\max_{x\in[0,L]}\{I^N(a)\}||U_0||^2  e^{2(1-\alpha)\gamma T} + \int_{0}^T\bigl(A_0^1 g_L\bigr)^2 dt}$\\\\\hline \\ \hline \\
		PDE:  & Non-conservative advection equation \\\\ \hline \\
		Equation & $u_t + au_x = 0,~~~ a(x)>0 ~;~~~ f_x=a(x)u_x$ \\\\
		Split version & $u_t + \alpha \bigl((au)_x-a_x u) + (1-\alpha)au_x = 0$\\\\
		& \hspace*{0.5cm}$\bullet$ $\alpha=0$: Non-conservative DG \\\\
		& \hspace*{0.5cm}$\bullet$ $\alpha=1/2$: Skew-symmetric DG \\\\
		& \hspace*{0.5cm}$\bullet$ $\alpha=1$: Conservative DG \\\\				  
		Continuous bound & $\displaystyle{\frac{||u(T)||^2}{\displaystyle{\max_{x\in[0,L]}\{a(x)\}}}  \leq \frac{||u_0||^2}{\displaystyle{\min_{x\in[0,L]}\{a(x)\}}} + \int_{0}^T g_L^2 dt}$\\\\
		Discrete bound  &$\displaystyle{\frac{||U(T)||^2}{\max_{x\in[0,L]}\{I^N(a)\}} \leq \frac{3||U_0||^2}{\min_{x\in[0,L]}\{I^N(a)\}} e^{2\alpha\gamma T} + \int_{0}^T g_L^2 dt}$\\\\\hline
	\end{tabu}
	\egroup
\end{table}

\section{Traditional energy bounds}\label{sec:Preliminaries}
The problem in \eqref{eq:general_advection_problem} is said to be \emph{strongly well-posed}, if it is well-posed and  the solution energy measured with the $L^2$ norm, $||u|| = \sqrt{\int_{0}^{L}u^{2}dx}$, satisfies the bound (see \cite{2014:Nordstrom}):

\begin{equation}
||u(\cdot,T)|| \leq K(T) \biggl(||u_0||^2 + \int_0^T \bigl(|g_L(\tau)|^2 + |g_R(\tau)|^2 \bigr)d\tau \biggr),
\end{equation}
where $K(t)$ is a constant that does not depend on the spatial coordinate, $x$. Since the source term $F(u,x)$ depends on the solution, its contribution on the integral $\int_0^T ||F||^2 d\tau$ is not included, although this source term will contribute to shape the energy amplitude, $K(t)$. For the one-dimensional advection equation, in which $f = a(x)u$, $F=0$ for its conservative version, and $F=a_xu$ for its non-conservative version, the precise expression for this bound was shown in \cite{Lorenz:1989fk,Kopriva2}:

\begin{equation}
||u(\cdot,T)||^2 + \int_0^T \bigl( \beta u^2(0,t) + \delta u^2(L,t)\bigr) d\tau  \leq e^{2\gamma T} \biggl( ||u_0||^2 + \int_0^T \bigl(\delta |g_L|^2 + \beta |g_R|^2 \bigr) d\tau \biggr),
\label{eq:traditional_bound}
\end{equation}
where the coefficients $\delta$ and $\beta$ depends on the behaviour of the left and right boundaries respectively. Hence, $(\delta,\beta)=1$ if they act as inflow, whereas $(\delta,\beta)=0$ when they act as outflow. The value of these two coefficients depend on the sign of the advection speed $a(x)$ at the boundaries. Additionally, note that the value of $\gamma$ depends on the bounds of the advection speed derivative:

\begin{equation}
\begin{split}
2\gamma = \min_{x \in [0,T]} \left| a_x\right|, &\text{ for the conservative form,}\\
2\gamma = \max_{x \in [0,T]} \left| a_x\right|, &\text{ for the non-conservative form.}
\end{split}	
\end{equation}

In this text we will consider the PDE when $a(x)$ does not change its sign (for the sake of simplicity, we will just consider $a(x)>0$). Nevertheless, the advection speed derivative $a_x$ is allowed to change sign within the computational domain, and $\gamma$ will be positive, allowing for an exponential growth of the energy. Whether this outcome is consistent or not with the analytical solution must be determined by means of alternative energy bounds. 


\section{Alternative continuous bounds}\label{sec:AlternativeContinuousBounds}
In this section, alternative energy bounds that show no exponential growth will actually be seen in the continuous setting. To do so, we define the following alternative inner products and norms:

\begin{equation}
\left\langle u,v\right\rangle_{a^\delta} = \int_{0}^T a(x)^{\delta} u(x)v(x) dx, ~~||u||_{a^\delta}^2 = \left\langle u,u\right\rangle_{a^\delta} \int_{0}a(x)^{\delta}u(x)^2 dx,
\label{eq:a_norm}
\end{equation}
where $\delta$ selects the $a-$norm ($\delta=1$), or the $1/a-$norm ($\delta=-1$). The $1/a-$norm has been previously used \cite{fornberg1975fourier,gottlieb2001spectral} to prove stability of the Fourier collocation method. 
Since $a(x)$ refers to the advection speed, which has been restricted to be positive, \eqref{eq:a_norm} satisfies all formal definitions of a norm. Both $a-$norm and $1/a-$norm can be related to the $L^2$ norm by

\begin{equation}
\min_{x\in[0,L]}\{a(x)\}||u||^2 \leq ||u||_a^2 \leq \max_{x\in[0,L]}\{a(x)\} ||u||^2,
\label{eq:L2_a^2norm}
\end{equation}
and,

\begin{equation}
\frac{||u||^2}{\displaystyle{\max_{x\in[0,L]}\{a(x)\}}} \leq ||u||_\frac{1}{a}^2 \leq \frac{||u||^2}{\displaystyle{\min_{x\in[0,L]}\{a(x)\}}},
\label{eq:L2_1/a^2norm}
\end{equation}
respectively. Thus, any bound measured with either the $a-$norm or the $1/a-$norm can be translated later to the $L^2$ norm. We will first obtain the continuous bound for the conservative advection equation,
\begin{equation}
u_t + [a(x)u]_x = 0,
\label{eq:conservative_eq}
\end{equation}
and that of the non-conservative equation,
\begin{equation}
u_t + a(x)u_x = 0
\label{eq:nonconservative_eq}
\end{equation}
will be derived afterwards. The domain, initial conditions and boundary conditions are those presented in \eqref{eq:general_advection_problem}.

\subsection{Continuous bound of the conservative advection equation}
We will first obtain the energy bound measured with the $a-$norm defined in \eqref{eq:a_norm}. To do so, instead of multiplying by $u(x,t)$ to obtain the energy estimate (i.e. the traditional form to perform this estimation), \eqref{eq:conservative_eq} is multiplied by the flux $f=a(x)u$:

\begin{equation}
\begin{split}
a(x)u u_t + f f_x&=0,\\
f(0,t) &= a(0)g_L.
\end{split}
\label{eq:conservative_flux}	
\end{equation}
To obtain the energy, \eqref{eq:conservative_flux} is integrated over the physical domain,

\begin{equation}
\int_{0}^L a(x)uu_t dx + \int_{0}^L f f_x dx = 0.
\label{eq:eqnIntegration}
\end{equation}
The first term in \eqref{eq:eqnIntegration} consists of the solution energy measured with the \textit{a-}norm, whilst the second term represents the contribution of the physical boundaries to the energy,

\begin{equation}
\int_{0}^L a(x) u u_t dx = \frac{1}{2}\frac{d}{dt}\int_{0}^L a(x) u^2dx = \frac{1}{2}\frac{d}{dt}||u||^2_{a} = -\frac{1}{2}f^2\biggr|_{0}^{L} = -\frac{1}{2}f_R^2 + \frac{1}{2}a(0)^2 g_L^2.
\end{equation}

The flux at the boundaries has been imposed according to the direction of the travelling physical waves. Thus, the solution energy satisfies
\begin{equation}
\frac{d}{dt}||u||_{a}^2+ f_R^2 = a(0)^2 g_L^2,
\end{equation}
and can be bounded by
\begin{equation}
\frac{d}{dt}||u||_{a}^2 \leq a(0)^2 g_L^2.
\label{eq:conservative_continuous_bound}
\end{equation}
Furthermore, time integration of \eqref{eq:conservative_continuous_bound} yields

\begin{equation}
||u(T)||_{a}^2 \leq ||u_0||^2_{a} + \int_{0}^T a(0)^2 g_L^2 dt.
\end{equation}
Using the relationship within the $L^2$ and the $a-$norm, \eqref{eq:L2_a^2norm}, we obtain the $L^2-$measured energy bound:

\begin{equation}
\min_{x\in[0,L]}\{a(x)\}||u(T)||^2 \leq \max_{x\in[0,L]}\{a(x)\}||u_0||^2 + \int_{0}^T a(0)^2 g_L^2 dt.\label{eq:conservative_continuous_bound_time}
\end{equation}
Thus, no exponential growth is expected in the analytical solution. The numerical scheme should mimic a discrete version this behaviour.

\subsection{Continuous bound of the non-conservative advection equation}
Regarding the non-conservative generic advection problem, its alternative energy estimate is obtained multiplying \eqref{eq:nonconservative_eq} by $u$, and then dividing by $a(x)$

\begin{equation}
\begin{split}
\frac{u u_t}{a(x)} + uu_x&=0,\\
u(0,t) &= g_L,
\end{split}
\label{eq:nonconservative_a}	
\end{equation}
which is nonsingular since $a(x)>0$. Following \eqref{eq:conservative_flux}, \eqref{eq:nonconservative_a} is integrated over the physical domain
\begin{equation}
\int_{0}^L \frac{uu_t}{a(x)}dx + \int_{0}^Luu_x dx = 0. \label{eq:nonconservative_energy_estimate_anorm}
\end{equation}
As in the conservative case, the first term is the solution energy, measured with the $1/a$-norm defined in \eqref{eq:a_norm} (with $\delta=-1$): 

\begin{equation}
\frac{1}{2}\frac{d}{dt}\int_{0}^L\frac{u^2}{a(x)}dx =\frac{1}{2}\frac{d}{dt} ||u||_{\frac{1}{a}}^2,
\end{equation}
The second term of \eqref{eq:nonconservative_energy_estimate_anorm} can be integrated by parts to show the energy entering and leaving the domain through the physical boundaries

\begin{equation}
\frac{1}{2}u^2\biggr|_{0}^L = -\frac{1}{2}g_L^2 + \frac{1}{2}u(L,t)^2.
\end{equation}
Putting it all together, the energy in terms of the $1/a-$norm varies as
\begin{equation}
\frac{d}{dt} ||u||_{\frac{1}{a}}^2+ u(L,t)^2 = g_L^2,
\end{equation}
which when time integrated leads to the bound
\begin{equation}
||u(T)||_{\frac{1}{a}}^2 \leq ||u_0||_{\frac{1}{a}}^2 + \int_0^T g_L^2dt.
\end{equation}
Finally, using the relationship between the $L^2$ and the $1/a-$norm, \eqref{eq:L2_1/a^2norm}, we get the energy bound
\begin{equation}
\frac{||u(T)||^2}{\displaystyle{\max_{x\in[0,L]}\{a(x)\}}}  \leq \frac{||u_0||^2}{\displaystyle{\min_{x\in[0,L]}\{a(x)\}}} + \int_{0}^T g_L^2 dt,
\label{eq:nonconservative_bound}
\end{equation}
which proves that the energy at any time is bounded, and does not exhibit any exponential growth. As before, the numerical scheme should mimic a discrete version of this behaviour.

\subsection{Summary}
Whilst traditional $L^2$ energy bounds show energy growth in general, $\gamma > 0$, the more strict $a-$norm estimate shows no growth when $a(x)>0$ (or equivalently $a(x)<0$) in both conservative and non-conservative forms of the equation. This result encourages us to use this approach to review the energy bounds of the discontinuous Galerkin approximation. Deriving the discrete version of the bounds, and bearing in mind that the final goal is to derive DG schemes that follow the continuous bound, gives us insight into the effect of the numerical errors on stability. We note in passing that alternative methods to obtain the energy bounds exist. The von Neumann analysis for non-constant advection speeds developed in \cite{2016:Manzanero} is also a powerful tool to understand the behaviour of numerical schemes.

\section{Discrete DGSEM-GL version of the bounds}\label{sec:DGSEMversion}
In this work, we consider the discontinuous Galerkin formulation with Gauss-Lobatto points. The scheme satisfies all the formal definitions of a summation-by-parts (SBP) simultaneous-approximation-term (SAT) scheme. This allows the removal of aliasing errors by means of split formulations (see \cite{2013:Gassner,2013:gass,2015:kopriva,2016:gassner}), and therefore, to obtain schemes that are provably stable. 

We will first introduce the notation used throughout this work. The physical domain $\Omega=\{x ~|~ x\in[0,L]\}$ is divided in $K$ non-overlapping elements $\Omega=\{x ~|~ x\in[x^k,x^{k+1}]\}$, in which the solution is approximated by $N$ degree polynomials (they will be said to belong to the $\mathcal{P}^N$ space). The interpolation is performed at the Gauss-Lobatto points $\{\xi_j\}_{j=0}^N$, which are then mapped individually from the local coordinate $\xi \in[-1,1]$ frame to each element domain by means of a linear one-dimensional mapping, $X^{el}(\xi)$

\begin{equation}
x_j^{el} = X^{el}(\xi_j) = x^{el} + \frac{1}{2}(x^{el+1}-x^{el})(\xi_j + 1)	.
\end{equation}
We adopt the system used in in \cite{Kopriva2} where capital symbols refer to the interpolated version of each variable. For instance, for the solution,
\begin{equation}
U^{el}(\xi) = I^N[u(X^{el}(\xi))]	 = \sum_{j=0}^N U_j^{el} l_j(\xi), ~~~U_j^{el} = u(x_j^{el}),
\end{equation}
and for the fluxes:

\begin{equation}
F^{el}(\xi) = I^N[f(X^{el}(\xi))]	 = \sum_{j=0}^N F_j^{el} l_j(\xi), ~~~F_j^{el} = f(x_j^{el}) = a(x_j^{el}) U_j^{el}.
\end{equation}
We will also obtain the discrete version of the advection speed, $a(x)$, $A^{el}(\xi)$ as

\begin{equation}
A^{el}(\xi) = I^N[a(X^{el})]	 = \sum_{j=0}^N A_j^{el} l_j(\xi), ~~~A_j^{el} = a(x_j^{el}).
\end{equation}
Note that because the Gauss-Lobatto points include the endpoints, the discrete version of the advection speed, $A^{el}(\xi)$, will remain continuous across the inter-element interfaces (assuming that $a(x)$ is continuous). 

The polynomial space is spanned by the Lagrange polynomials, $l_j(\xi)$, with nodes at the Gauss-Lobatto points
\begin{equation}
l_j(\xi) = \prod_{\substack{m=0 \\ m \neq j}}^{N} \frac{\xi - \xi_m}{\xi_j - \xi_m}.
\end{equation}
This basis allows one to compute the required derivatives by means of a derivative matrix, $[\boldsymbol{D}]$, defined as
\begin{equation}
D_{ij} = l'_j(\xi_i).
\end{equation}

We also adopt the following matrix-vector form: A vector contains the nodal degrees of freedom of a certain variable, for instance $\{\underline{\boldsymbol{A}}^{el}\}$. The notation $[\boldsymbol{A}^{el}]=\text{diag}(\{\underline{\boldsymbol{A}}^{el}\})$) represents the diagonal matrix whose diagonal entries are those nodal values. This allows us to compactly write the following terms, usually arising from the split formulation, as
\begin{equation}
\begin{split}
I^N[A^{el}U^{el}_\xi] &=\{\underline{\boldsymbol{l}}^{el}(\xi)\}^T[\boldsymbol{A}^{el}][\boldsymbol{D}^{el}] \{\underline{\boldsymbol{U}}^{el}\},\\
I^N[A^{el}_\xi U^{el}]&=\{\underline{\boldsymbol{l}}^{el}(\xi)\}^T [\boldsymbol{A}^{el}_\xi] \{\underline{\boldsymbol{U}}^{el}\},~~~[\boldsymbol{A}^{el}_\xi] = \text{diag}([\boldsymbol{D}^{el}]\{\underline{\boldsymbol{A}}^{el}\})\\
(I^N[A^{el}U^{el}])_\xi &=\{\underline{\boldsymbol{l}}^{el}(\xi)\}^T[\boldsymbol{D}^{el}][\boldsymbol{A}^{el}] \{\underline{\boldsymbol{U}}^{el}\}.\\
\end{split}
\label{eq:matricial_forms}
\end{equation}

We will discretise the advection equation, either in conservative or non-conservative form, by means of a general split formulation for the flux, $f=au$. The equation solved depends a the parameter $\theta$, namely
\begin{equation}
u_t + f_x = \theta a_x u	.
\label{eq:EqnwithTheta}
\end{equation}
Setting $\theta=0$ solves the conservative equation, while the case $\theta=1$ recovers the non-conservative equation. The split form is described in \cite{Kopriva2},
\begin{equation}
u_t + \alpha f_x + (1-\alpha)(a_x u + au_x) = \theta a_x u.
\label{eq:skew_symmetric_advection}
\end{equation}
Equation \eqref{eq:skew_symmetric_advection} is multiplied by a test function, $\Phi\in\mathcal{P}^N$, and integrated in each element
to get the weak form
\begin{equation}
\begin{split}
\frac{\Delta x^{el}}{2}\left\langle\Phi^{el},u_t\right\rangle + \alpha \left\langle\Phi^{el},f_\xi\right\rangle &\\
+ (1-\alpha)\bigl(\left\langle\Phi^{el},a_\xi u\right\rangle + \left\langle au_\xi,\Phi^{el}\right\rangle\bigr) &= \theta \left\langle\Phi^{el},a_\xi u\right\rangle,
\end{split}
\label{eq:FirstWeakForm}
\end{equation}
where the inner product is $\left\langle u,v\right\rangle=\int_{-1}^{1}uvd\xi$.

The second and fourth terms of \eqref{eq:FirstWeakForm} are integrated by parts, and the interface fluxes that appear are replaced by a numerical flux $f^\star$:
\begin{equation}
\begin{split}
\frac{\Delta x^{el}}{2}\left\langle\Phi^{el},u_t\right\rangle + f^\star \Phi^{el}\biggr|_{-1}^1 - \alpha \left\langle\Phi^{el}_\xi,f\right\rangle &\\
+ (1-\alpha)\bigl(\left\langle\Phi^{el},a_\xi u\right\rangle - \left\langle u,(a\Phi^{el})_\xi\right\rangle\bigr) &= \theta \left\langle\Phi^{el},a_\xi u\right\rangle.
\end{split}
\end{equation}
Lastly, inner product integrals are computed with numerical quadratures (i.e. Gauss-Lobatto points). Thus, we replace the inner products $\left\langle\cdot,\cdot\right\rangle$ by their numerical version $\left\langle\cdot,\cdot\right\rangle_N$, and the arguments by their polynomial approximations
\begin{equation}
\begin{split}
\frac{\Delta x^{el}}{2}\left\langle\Phi^{el},U^{el}_t\right\rangle_N + F^\star \Phi^{el}\biggr|_{-1}^1 - \alpha \left\langle\Phi^{el}_\xi,F^{el}\right\rangle_N& \\
+ (1-\alpha)\bigl(\left\langle\Phi^{el},I^N[A^{el}_\xi U^{el}]\right\rangle_N - \left\langle U^{el},(I^N[A^{el}\Phi^{el}])_\xi\right\rangle_N\bigr) &= \theta \left\langle\Phi^{el},I^N[A^{el}_\xi U^{el}]\right\rangle_N.
\end{split}
\label{eq:discrete_dg_scheme}
\end{equation}

Using $N+1$ linearly independent test functions (e.g. the Lagrange polynomials), we obtain the differential equations for the solution discrete degrees of freedom, $U_j^{el}(t)$. Henceforth, for the sake of simplicity, we will drop the $el-$ index. 

The discrete formulation is selected with the split form coefficient, $\alpha$, where $\alpha=0$ gives the non-conservative DG, $\alpha=1/2$ gives to the skew-symmetric DG, and $\alpha=1$ gives the conservative DG. It is also possible to switch between the two forms of the equation through the parameter $\theta$, where with $\theta=0$ we solve the conservative equation, whilst with $\theta=1$ we recover the non-conservative equation. 

Now we proceed to obtain the discrete version of the continuous bounds in \eqref{eq:conservative_continuous_bound_time} and \eqref{eq:nonconservative_bound}. To do so, we define the discrete version of the $a-$norm, the $a-$inner product, the $1/a$-norm, and the $1/a-$inner product as
\begin{equation}
\begin{split}
\left\langle U,V\right\rangle_{a^\delta,N} = \sum_{m=0}^N w_m A_m^\delta U_m V_m =\{\underline{\boldsymbol{U}}\}^T[\boldsymbol{M}][\boldsymbol{A}^\delta]\{\underline{\boldsymbol{V}}\}, \text{ and}\\
||U||_{a^\delta}^2 = \left\langle U,U\right\rangle_{a^\delta,N} = \sum_{m=0}^N w_m A_m^\delta U_m^2=\{\underline{\boldsymbol{U}}\}^T[\boldsymbol{M}][\boldsymbol{A}^\delta]\{\underline{\boldsymbol{U}}\},\\
\end{split}
\label{eq:discrete_a_norms}
\end{equation}
where $\delta=1$ refers to the $a-$norm, and $\delta = -1$ refers to the $1/a-$norm. In \eqref{eq:discrete_a_norms}, the $w_m$ are the Gauss-Lobatto quadrature weights, and $[\boldsymbol{M}]$ the mass matrix, whose entries are the quadrature weights placed along the main diagonal. The two discrete norms can be also related to the continuous $L^2$ norm. First, the discrete $a^\delta-$ norm can be related to the Gauss-Lobatto discrete norm by
\begin{equation}
\min\{A_j\}_{j=0}^N ||U||_N^2\leq ||U||_{a,N}^2\leq\max\{A_j\}_{j=0}^N ||U||_N^2,
\end{equation}
and related to the $L^2$ norm (see Section 5.3 in \cite{Canuto:2006}) as
\begin{equation}
\min\{A_j\}_{j=0}^N ||U||^2\leq ||U||_{a,N}^2\leq 3\max\{A_j\}_{j=0}^N ||U||^2.
\label{eq:continuous_discrete_anorm}
\end{equation}
equation \eqref{eq:continuous_discrete_anorm} is valid as long as $A^{el}(\xi)> 0$ in each element. Similarly, for the $1/a-$norm:
\begin{equation}
\frac{||U||^2}{\max\{A_j\}_{j=0}^N} \leq ||U||_{\frac{1}{a},N}^2\leq \frac{3||U||^2}{\min\{A_j\}_{j=0}^N} .
\end{equation}
Finally, the $a^\delta$ norm and inner products are extended to the whole domain by summing the all the elemental contributions
\begin{equation}
\begin{split}
\left\langle U,V\right\rangle_{a^\delta,N} &=\frac{\Delta x}{2}\sum_{el=1}^K \left\langle U^{el},V^{el}\right\rangle_{a^\delta,N},\\
||U||_{a^\delta}^2 &=\frac{\Delta x}{2}\sum_{el=1}^K ||U^{el}||_{a^\delta}^2.
\end{split}
\end{equation}

\subsection{Discrete bound of the conservative advection equation}
We set $\theta=0$ to recover the conservative equation and rearrange \eqref{eq:discrete_dg_scheme} as the conservative or standard DG (see \cite{2009:Kopriva}), plus a correction term that arises from the split formulation
\begin{equation}
\begin{split}
&\frac{\Delta x}{2}\left\langle\Phi,U_t\right\rangle_N + F^\star \Phi\biggr|_{-1}^1 - \left\langle\Phi_\xi,F\right\rangle_N \\
+& (1-\alpha)\bigl(\left\langle\Phi_\xi,F\right\rangle_N + \left\langle\Phi,I^N[A_\xi U]\right\rangle_N - \left\langle U,(I^N[A\Phi])_\xi\right\rangle_N\bigr) =0.
\end{split}
\label{eq:conservative_discrete_dg_scheme}
\end{equation}
Hence, the precise form of the correction term is
\begin{equation}
\left\langle\Phi_\xi,F\right\rangle _N + \left\langle \Phi,I^N[A_\xi U]\right\rangle _N - \left\langle U,(I^N[A\Phi])_\xi\right\rangle _N = N_1 + N_2 + N_3.
\label{eq:skew_correction}
\end{equation}

To derive an energy bound, similar to its continuous counterpart, the test function is replaced by the discrete version of the fluxes, $\Phi=F$
\begin{equation}
\begin{split}
&\frac{\Delta x}{2}\left\langle F,U_t\right\rangle _N + F^\star F\biggr|_{-1}^1 - \left\langle F_\xi,F\right\rangle _N \\
+& (1-\alpha)\bigl(\left\langle F_\xi,F\right\rangle _N + \left\langle F,I^N[A_\xi U]\right\rangle _N - \left\langle U,(I^N[AF])_\xi\right\rangle _N\bigr) =0.
\end{split}
\label{eq:conservative_discrete_bound_1}
\end{equation}
The first term in \eqref{eq:conservative_discrete_bound_1} reproduces the time derivative of the discrete energy measured with the $a-$norm,
\begin{equation}
\left\langle F,U_t\right\rangle _N = \left\langle AU,U_t\right\rangle _N = \left\langle U,U_t\right\rangle _{a,N} = \frac{1}{2}\frac{d}{dt}||U||^2_{a,N}.
\label{eq:conservative_energy_timederivative}
\end{equation}
Regarding the third term in \eqref{eq:conservative_discrete_bound_1}, the summation-by-parts property (see \cite{2013:Carpenter}) holds, and thus
\begin{equation}
\left\langle F_\xi,F\right\rangle _N = \frac{1}{2}\bigl(F\bigr)^2\biggl|_{-1}^1.
\end{equation}
Lastly, in the correction term defined in \eqref{eq:skew_correction}, both first and third term can be rewritten (following the summation-by-parts property) as
\begin{equation}
\begin{split}
N_1 + N_3 &= \left\langle F_\xi,F\right\rangle _N  - \left\langle U,(I^N[AF])_\xi\right\rangle _N\\
&= 	- \left\langle F_\xi,F\right\rangle _N + \left\langle U_\xi,I^N[AF]\right\rangle _N,
\end{split}
\end{equation}
in which boundary terms arising from the summation-by-parts cancel so only volume integrals contribute to the estimate. All three terms involved in the correction term can be converted to the $a-$inner product. The first is
\begin{equation}
- \left\langle F_\xi,F\right\rangle _N = -\left\langle F_\xi , U\right\rangle _{a,N} = -\left\langle \bigl(I^N[AU]\bigr)_\xi , U\right\rangle _{a,N}.
\end{equation}
The second
is
\begin{equation}
\left\langle F,I^N[A_\xi U]\right\rangle _N = 	\left\langle U,I^N[A_\xi U]\right\rangle _{a,N}.
\end{equation}
And the third becomes
\begin{equation}
\left\langle U_\xi,I^N[AF]\right\rangle _N = \left\langle I^N[A U_\xi],U]\right\rangle _{a,N}.
\end{equation}
Therefore, the correction term, $\mathcal{L}(A,U)$, is
\begin{equation}
\mathcal{L}(A,U) = \left\langle I^N[A_\xi U] + I^N[AU_\xi] - \bigl(I^N[AU]\bigr)_\xi , U\right\rangle _{a,N},
\end{equation}
which does not vanish since the product derivative rule does not have a discrete equivalent. This term represents the aliasing errors introduced in the discrete weak formulation of the original equation, which are projected (with the a-norm inner product) onto the solution. We can bound these errors, since, using the matrix form shown in \eqref{eq:matricial_forms} we can write the inner product as
\begin{equation}
\mathcal{L}(A,U) = \{\underline{\boldsymbol{U}}\}^T[\boldsymbol{A}][\boldsymbol{M}]\bigl([\boldsymbol{A}_\xi] + [\boldsymbol{A}][\boldsymbol{D}]-[\boldsymbol{D}][\boldsymbol{A}] \bigr)\{\underline{\boldsymbol{U}}\}.
\end{equation}
Thus, we can use the Cauchy-Schwartz inequality using the $a-$norm to bound the aliasing term as
\begin{equation}
\mathcal{L}(A,U) \le  \gamma ||U||_{a,N}^2,
\label{eq:aliasing_bound_cons}
\end{equation}
where
\begin{equation}
\gamma = \bigl\Vert[\boldsymbol{A}][\boldsymbol{M}]\bigl([\boldsymbol{A}_\xi] + [\boldsymbol{A}][\boldsymbol{D}]-[\boldsymbol{D}][\boldsymbol{A}]\bigr)\bigr\Vert_{a,N},
\label{eq:gamma_definition}
\end{equation}
The coefficient $\gamma$ depends on the advection speed, $A(\xi)$, and the polynomial degree, $N$.  Gathering all the terms together, we get the elemental contribution to the time derivative of the energy,
\begin{equation}
\frac{\Delta x}{2}\frac{d}{dt}||U^{el}||^2_{a,N} + \bigl(2F^\star-F^{el}\bigr) F^{el}\biggr|_{-1}^1 \leq 2(1-\alpha)\gamma^{el}||U^{el}||_{a,N}^2.
\label{eq:conservative_discrete_bound_2}
\end{equation}
Note that now we specifically include the $el-$ index. Next, \eqref{eq:conservative_discrete_bound_2} is summed over all elements, to obtain the total energy,
\begin{equation}
\frac{d}{dt}||U||^2_{a,N} + \sum_{el=1}^K\bigl(2F^\star-F^{el}\bigr) F^{el}\biggr|_{-1}^1 \leq 2(1-\alpha)\gamma ||U||^2_{a,N},
\end{equation}
where $\gamma$ is bounded by the largest value over all the elements,
\begin{equation}
\gamma = \max_{ el} \frac{\gamma^{el}}{\Delta x^{el}/2}.
\end{equation}

The interior interfaces contribution to the total energy vanishes as long as central fluxes are considered (see \cite{Kopriva2}), wheras the physical boundary contributions, computed with upwind fluxes, are
\begin{equation}
\sum_{el=1}^K\bigl(2F^\star-F^{el}\bigr) F^{el}\biggr|_{-1}^1 = \bigl(F^K(1)\bigr)^2 - \bigl(2A_0^1g_L -F^1(-1)\bigr)F^1(-1).
\label{eq:cons_fluxes_1}
\end{equation}
The upwind flux stabilises the inflow condition, since
\begin{equation}
\bigl(2A_0^1g_L -F^1(-1)\bigr)F^1(-1) = \bigl(A_0^1 g_L\bigr)^2-\bigl(A_0^1 g_L - F^{1}(-1)\bigr)^2.
\label{eq:cons_fluxes_2}	
\end{equation}
Thus, we get the energy estimate
\begin{equation}
\frac{d}{dt}||U||^2_{a,N} +\bigl(A_0^1 g_L - F^{1}(-1)\bigr)^2 + \bigl(F^K(1)\bigr)^2  \le 2(1-\alpha)\gamma||U||_{a,N}^2+ \bigl(A_0^1 g_L\bigr)^2,
\label{eq:conservative_discrete_estimate}
\end{equation}
with the upper bound
\begin{equation}
\frac{d}{dt}||U||^2_{a,N} \leq  2(1-\alpha)\gamma||U||_{a,N}^2+ \bigl(A_0^1 g_L\bigr)^2.
\label{eq:conservative_discrete_bound}
\end{equation}
Time integration of \eqref{eq:conservative_discrete_bound} yields
\begin{equation}
||U(T)||^2_{a,N} \leq ||U_0||^2_{a,N}e^{2(1-\alpha)\gamma T} + \int_{0}^T\bigl(A_0^1 g_L\bigr)^2 dt,
\end{equation}
which is related to the continuous $L_2$ norm using \eqref{eq:continuous_discrete_anorm} as
\begin{equation}
\min\{A_j\}_{j=0}^N||U(T)||^2 \leq 3\max\{A_j\}_{j=0}^N||U_0||^2  e^{2(1-\alpha)\gamma T} + \int_{0}^T\bigl(A_0^1 g_L\bigr)^2 dt. \label{eq:conservative_L2_bound}
\end{equation}

As a conclusion, the final outcome is that aliasing errors drive the instability, causing an exponential growth. Choosing the correct value of the parameter $\alpha$ makes it possible to remove those errors from the energy estimate. Precisely, for the conservative equation, the parameter $\alpha$ should be equal to $1$, that is, a conservative DG scheme. If it is the case, has been proven that the numerical solution will not exhibit exponential growth, as does the physical solution.
\subsection{Discrete bound of the non-conservative advection equation}
We now switch $\theta$ in \eqref{eq:EqnwithTheta} to $1$ to get the non-conservative equation. As with \eqref{eq:conservative_discrete_dg_scheme}, we rearrange \eqref{eq:discrete_dg_scheme} to be regarded as the non-conservative or standard DG (see \cite{2009:Kopriva}), plus a correction term that arises from the split formulation,
\begin{equation}
\begin{split}
&\frac{\Delta x}{2}\left\langle \Phi,U_t\right\rangle _N + F^\star \Phi\biggr|_{-1}^1 - \left\langle \Phi_\xi,F\right\rangle _N = \left\langle \Phi,I^N[A_\xi U]\right\rangle _N  \\
-& (1-\alpha)\bigl(\left\langle \Phi_\xi,F\right\rangle _N + \left\langle \Phi,I^N[A_\xi U]\right\rangle _N - \left\langle U,(I^N[A\Phi])_\xi\right\rangle _N\bigr).
\end{split}
\label{eq:nonconservative_discrete_dg_scheme}
\end{equation}
In \eqref{eq:nonconservative_discrete_dg_scheme}, the correction term is defined as
\begin{equation}
\left\langle \Phi_\xi,F\right\rangle _N + \left\langle \Phi,I^N[A_\xi U]\right\rangle _N - \left\langle U,(I^N[A\Phi])_\xi\right\rangle _N = N_1 + N_2 + N_3.
\label{eq:correction_nonconservative}
\end{equation}

This time the discrete energy bound is derived by replacing the test function by the discrete quotient between the solution and the advection speed, $\Phi=I^N[U/A]$ to get
\begin{equation}
\begin{split}
&\frac{\Delta x}{2}\left\langle I^N\biggl[\frac{U}{A}\biggr],U_t\right\rangle _N + F^\star \biggl[\frac{U}{A}\biggr]\biggr|_{-1}^1 + \left\langle U,(I^N[AI^N\biggl[\frac{U}{A}\biggr]])_\xi\right\rangle _N =  \\
+& \alpha\bigl(\left\langle \biggl(I^N\biggl[\frac{U}{A}\biggr]\biggr)_\xi,F\right\rangle _N + \left\langle I^N\biggl[\frac{U}{A}\biggr],I^N[A_\xi U]\right\rangle _N - \left\langle U,(I^N[AI^N\biggl[\frac{U}{A}\biggr]])_\xi\right\rangle _N\bigr),
\end{split}
\label{eq:nonconservative_discrete_bound_1}
\end{equation}
where the correction term has been rearranged. The first term in \eqref{eq:nonconservative_discrete_bound_1} reproduces the time derivative of the discrete energy measured with the $1/a-$norm:
\begin{equation}
\left\langle I^N\biggl[\frac{U}{A}\biggr],U_t\right\rangle _N = \left\langle U,U_t\right\rangle _{\frac{1}{a},N} = \frac{1}{2}\frac{d}{dt}||U||^2_{\frac{1}{a},N}.
\end{equation}
We can use the summation-by-parts property for the third term in \eqref{eq:nonconservative_discrete_bound_1} to write the volume term in terms
of surface quantities
\begin{equation}
\left\langle U,(I^N[AI^N\biggl[\frac{U}{A}\biggr]])_\xi\right\rangle _N =  \left\langle U,U_\xi\right\rangle _N = \frac{1}{2}\bigl(U\bigr)^2\biggr|_{-1}^1.
\end{equation}

Both first and third terms of the correction term defined in \eqref{eq:correction_nonconservative} can be rewritten (following the summation-by-parts property) as
\begin{equation}
\begin{split}
N_1 + N_3 &= \left\langle \biggl(I^N\biggl[\frac{U}{A}\biggr]\biggr)_\xi,F\right\rangle _N - \left\langle U,(I^N[AI^N\biggl[\frac{U}{A}\biggr]])_\xi\right\rangle _N \\
&=  \left\langle I^N\biggl[\frac{U}{A}\biggr],AU_\xi\right\rangle _N - \left\langle I^N\biggl[\frac{U}{A}\biggr],F_\xi\right\rangle _N .
\end{split}
\end{equation}
Again, the three terms in the correction term can be written in terms of the $1/a-$inner product. We write the first as
\begin{equation}
- \left\langle I^N\biggl[\frac{U}{A}\biggr],F_\xi\right\rangle _N = - \left\langle U,F_\xi\right\rangle _{\frac{1}{a},N} ,
\end{equation}
the second as
\begin{equation}
\left\langle I^N\biggl[\frac{U}{A}\biggr],I^N[A_\xi U]\right\rangle _N = \left\langle U , I^N[A_\xi U]\right\rangle _{\frac{1}{a},N},
\end{equation}
and the third as
\begin{equation}
\left\langle I^N\biggl[\frac{U}{A}\biggr],AU_\xi\right\rangle _N = \left\langle U,AU_\xi\right\rangle _{\frac{1}{a},N}.
\end{equation}
Therefore, the correction term, $\mathcal{L}(A,U)$, is
\begin{equation}
\mathcal{L}(A,U) = \left\langle I^N[A_\xi U] + I^N[A U_\xi] - \bigl(I^N[A U]\bigr)_\xi , U\right\rangle _{\frac{1}{a},N},
\end{equation}
which again does not vanish since the product derivative rule does not have a discrete equivalent. This term represents the aliasing errors incurred in the discrete weak formulation of the original equation. Moreover, following a similar approach to \eqref{eq:aliasing_bound_cons}, it can be bounded as
\begin{equation}
{\mathcal{L}}(A,U) \leq  \gamma||U||_{\frac{1}{a},N}^2,
\end{equation}
where the definition of $\gamma$ is now based on the $1/a-$norm
\begin{equation}
\gamma = \bigl\Vert[\boldsymbol{A}]^{-1}[\boldsymbol{M}]\bigl([\boldsymbol{A}_\xi] + [\boldsymbol{A}][\boldsymbol{D}]-[\boldsymbol{D}][\boldsymbol{A}]\bigr)\bigr\Vert_{\frac{1}{a},N}.
\end{equation}
%
%
%
%

Putting it all together, the estimate of the time derivative of the energy inside each element (where the $el-$ index is again explicitly written) is 
\begin{equation}
\frac{\Delta x}{2}\frac{d}{dt}||U^{el}||^2_{\frac{1}{a},N} + \bigl(2U^\star-U^{el}\bigr) U^{el}\biggr|_{-1}^1 \le 2\alpha \gamma^{el}||U^{el}||_{\frac{1}{a},N}^2,
\label{eq:nonconservative_discrete_bound_2}
\end{equation}
where $U^\star$ is the interelement flux divided by the interelement advection speed. Next, \eqref{eq:nonconservative_discrete_bound_2} is summed across all elements, to obtain the time derivative of the total energy,
\begin{equation}
\frac{d}{dt}||U||^2_{\frac{1}{a},N} + \sum_{el=1}^K\bigl(2U^\star-U^{el}\bigr)U^{el}\biggr|_{-1}^1 \le 2\alpha \gamma^{el}||U^{el}||_{\frac{1}{a},N}^2.
\end{equation}
The contribution to the total energy of the interior interfaces vanishes when central fluxes are used (see \cite{Kopriva2}). When the physical boundary contributions are computed with upwind fluxes, we get terms like those in \eqref{eq:cons_fluxes_1} and \eqref{eq:cons_fluxes_2}. Thus,
\begin{equation}
\frac{d}{dt}||U||^2_{\frac{1}{a},N} +\bigl( g_L - U^{1}(-1)\bigr)^2 + \bigl(U^K(1)\bigr)^2 \le  2\alpha \gamma||U||_{\frac{1}{a},N}^2+ \bigl(g_L\bigr)^2,
\label{eq:nonconservative_discrete_estimate}
\end{equation}
which has an upper bound
\begin{equation}
\frac{d}{dt}||U||^2_{\frac{1}{a},N} \leq 2\alpha \gamma||U||_{\frac{1}{a},N}^2+ \bigl(g_L\bigr)^2.
\label{eq:nonconservative_discrete_bound}
\end{equation}
The energy at $t=T$ is obtained integrating \eqref{eq:nonconservative_discrete_bound} in time
\begin{equation}
||U(T)||^2_{\frac{1}{a},N} \leq ||U_0||^2_{\frac{1}{a},N}e^{2\alpha\gamma T} + \int_{0}^T g_L^2 dt,
\end{equation}
which can be related to the continuous $L_2$ norm using \eqref{eq:continuous_discrete_anorm},
\begin{equation}
\frac{||U(T)||^2}{\max\{A_j\}_{j=0}^N} \leq \frac{3||U_0||^2}{\min\{A_j\}_{j=0}^N} e^{2\alpha\gamma T} + \int_{0}^T g_L^2 dt.\label{eq:nonconservative_L2_bound}
\end{equation}

Therefore, we can conclude that, like in the conservative equation approximation, any exponential growth experienced by the numerical solution is due to aliasing errors in the flux interpolation. Notwithstanding, we can remove these errors by choosing the precise split form coefficient. In particular, the parameter $\alpha$ should be set to zero (i.e. a non-conservative discretisation, $\alpha=0$) to follow the behaviour of the analytical solution. 

\subsection{Summary of alternative discrete energy bounds}
We have obtained the discrete versions of the continuous bounds derived in Section \ref{sec:AlternativeContinuousBounds} for both conservative and non-conservative DG. These bounds have been derived using a split-form discontinuous Galerkin method with Gauss-Lobatto points. These alternative bounds make it possible to analyse the aliasing error removal by means of the split operator. These aliasing errors, which arise from the fact that the product rule does not hold, drive an exponential energy growth, which can be controlled with the parameter $\alpha$. Precisely, when solving the conservative equation, it suffices to select $\alpha=1$ (i.e. conservative DG), whilst selecting $\alpha=0$ (non-conservative DG) when solving the non-conservative equation is enough to prevent aliasing driven exponential growth.
It should be noticed that the alternative bounds derived in Section \ref{sec:AlternativeContinuousBounds} are only valid for advection speeds $a(x)$ that do not change sign in the domain. In the general case, strongly-stable schemes are recovered for $\alpha=1/2$ as shown in \cite{Kopriva2}.

\section{Effect of over-integration}\label{sec:overintegration}
Over-integration, also known as polynomial de-aliasing or consistent integration, is a technique commonly used to reduce aliasing in non-linear equations \cite{2003:Kirby}. In this work, we apply this methodology to the non-constant speed advection equation and find that, despite of being capable of removing the non-linear fluxes aliasing errors, it  introduces additional instabilities through the discrete weak form.

Since we are considering the non-constant speed advection equation, the discrete flux is the product of two $N-$degree polynomials, once the generally non-polynomial advection speed, $a(x)$, has been projected onto the solution space. Therefore, we may change the weak formulation to
\begin{equation}
\frac{\Delta x}{2}\left\langle \Phi,U_t\right\rangle _N + F^\star \Phi\biggr|_{-1}^1 - \left\langle \Phi_\xi,F\right\rangle _M = 0, \\
\label{eq:conservative_overInt}
\end{equation}
where $F$ is now a polynomial of degree $2N$ (since $F=AU$), and $M$ is the required number of quadrature points. For the Gauss-Lobatto points, $M> 3N/2$ is required to avoid inexact quadratures. Also, note that any of the split forms we have considered are algebraically identical to the conservative formulation, since all volume integrals are approximated and computed exactly. 

To compute the energy estimate, we replace $\Phi$ by $I^N[F]$ since the test function should belong to the solution space, $\mathcal{P}^N$,
\begin{equation}
\frac{\Delta x}{2}\left\langle I^N[F],U_t\right\rangle _N + F^\star F\biggr|_{-1}^1 - \left\langle \bigl(I^N[F]\bigr)_\xi,F\right\rangle _M = 0. \\
\label{eq:conservative_overInt}
\end{equation}
The interpolation operator $I^N[\bullet]$ has been omitted for the surface integral since Gauss-Lobatto points include the endpoints, and therefore both are identical (i.e. interpolation is not necessary at the boundaries). Following \eqref{eq:conservative_energy_timederivative}, we write the first term in \eqref{eq:conservative_overInt} as the time derivative of the energy, measured with the $a-$norm, whereas the last term is summated-by-parts
\begin{equation}
\begin{split}
\left\langle \bigl(I^N[F]\bigr)_\xi,F\right\rangle _M &= \left\langle F_\xi,F\right\rangle _M + \left\langle \bigl(I^N[F]\bigr)_\xi-F_\xi,F\right\rangle _M \\
&=\frac{1}{2}F\biggr|_{-1}^{+1} + \left\langle \bigl(I^N[F]\bigr)_\xi-F_\xi,F\right\rangle _M.
\end{split}
\end{equation}
Thus, the energy bound reads

\begin{equation}
\frac{\Delta x}{2}\frac{d}{dt}||U||^2_{a,N} + \bigl(2F^\star-F\bigr) F\biggr|_{-1}^1  = \left\langle \bigl(F - I^N[F]\bigr)_\xi,F\right\rangle _M.
\end{equation}
Which summed over all mesh elements gives
\begin{equation}
\frac{d}{dt}||U||^2_{a,N}  \leq (A_0^N g_L)^2 + \left\langle \bigl(F - I^N[F]\bigr)_\xi,F\right\rangle _M.
\label{eq:bound_overintegration}
\end{equation}
The last term on the right hand side represents the aliasing errors incurred when using the interpolant of the flux as the test function. Therefore, computing  the integrals exactly may not stabilise the scheme in this norm. 

\section{Numerical experiments}\label{sec:numericalExperiments}

We will show the accuracy of the bounds in \eqref{eq:conservative_L2_bound}, \eqref{eq:nonconservative_L2_bound}, and \eqref{eq:bound_overintegration} by examining the eigenvalues of the numerical implementations of the advection equation. This particular example considers the domain $x\in[-1,1]$ in which the advection speed is
\begin{equation}
a(x) = 1+(1-x^2)^5, \label{eq:a(x)}
\end{equation}
and with periodic boundary conditions at the two endpoints. This example (which comes from Hesthaven and Warbuton's book \cite{2008:hest}) is interesting since it contains large spectral content. 

We use the split discontinuous Galerkin method shown in \eqref{eq:discrete_dg_scheme} to solve both conservative and non-conservative equations. Moreover, we use the central fluxes,
\begin{equation}
F^\star(u_L,u_R) = \frac{a_L u_L + a_R u_R}{2},
\end{equation}
for both interior and physical boundaries, since we have enforced periodic boundary conditions. As a consequence of performing the interpolation with Gauss-Lobatto points, the discrete version of $a(x)$ will be continuous, and $a_L=a_R$, equal to the advection speed at the boundary. 

Replacing the test function $\Phi^{el}$ by the set of Lagrange polynomials $\{l_j\}_{j=0}^N$ in \eqref{eq:discrete_dg_scheme} and computing the inner products, one arrives to the following matrix system for each element
\begin{equation}
\frac{\Delta x^{el}}{2}[\boldsymbol{M}]\frac{d\{\underline{\boldsymbol{U}}^{el}\}}{dt} = [\boldsymbol{L}^{el}]\{\underline{\boldsymbol{U}}^{el-1}\} + [\boldsymbol{C}^{el}]\{\underline{\boldsymbol{U}}^{el}\} + [\boldsymbol{R}^{el}]\{\underline{\boldsymbol{U}}^{el+1}\},
\end{equation}
where the matrices $[\boldsymbol{L}^{el}]$, $[\boldsymbol{C}^{el}]$, and $[\boldsymbol{R}^{el}]$ are constant, i.e. they do not depend on $\{\boldsymbol{\underline{U}}\}$, and they just depend on the polynomial order and the discrete advection speed in each element, $A^{el}(\xi)$. Their precise expression is
\begin{align}
[\boldsymbol{L}^{el}] &= \frac{a^{el-1}}{2}\{\boldsymbol{\underline{l}}(-1)\}\{\boldsymbol{\underline{l}}(1)\}^T,\nonumber\\
[\boldsymbol{C}^{el}] &= \frac{a^{el-1}}{2}\{\boldsymbol{\underline{l}}(1)\}\{\boldsymbol{\underline{l}}(1)\}^T - \frac{a^{el}}{2}\{\boldsymbol{\underline{l}}(-1)\}\{\boldsymbol{\underline{l}}(-1)\}^T\label{LCRmatrices}\\
&+ \alpha[\boldsymbol{D}]^T[\boldsymbol{M}][\boldsymbol{A}^{el}] -(1-\theta-\alpha) [\boldsymbol{A}_\xi^{el}] + (1-\alpha) [\boldsymbol{A}^{el}][\boldsymbol{D}]^T[\boldsymbol{M}],\nonumber\\
[\boldsymbol{R}^{el}] &=- \frac{a^{el} }{2}\{\boldsymbol{\underline{l}}(1)\}\{\boldsymbol{\underline{l}}(-1)\}^T,\nonumber
\end{align}
where $a^{el-1}$ is the advection speed to the left of the element boundary, and $a^{el}$ is that to the right. Recall that $\theta=0$ generates the conservative equation, whilst $\theta=1$ generates the non-conservative equation.\\

The stability of the scheme is studied using the system of eigenvalues over all elements. We build a mesh with $K=200$ elements and polynomial order $N=5$. We have considered the three relevant cases for the split operator coefficient: conservative DG ($\alpha=1$), skew-symmetric DG ($\alpha=1/2$), and non-conservative DG ($\alpha=0$). 

\begin{figure}
	\centering
	\subfigure[Conservative equation spectra]{\includegraphics[scale=0.23]{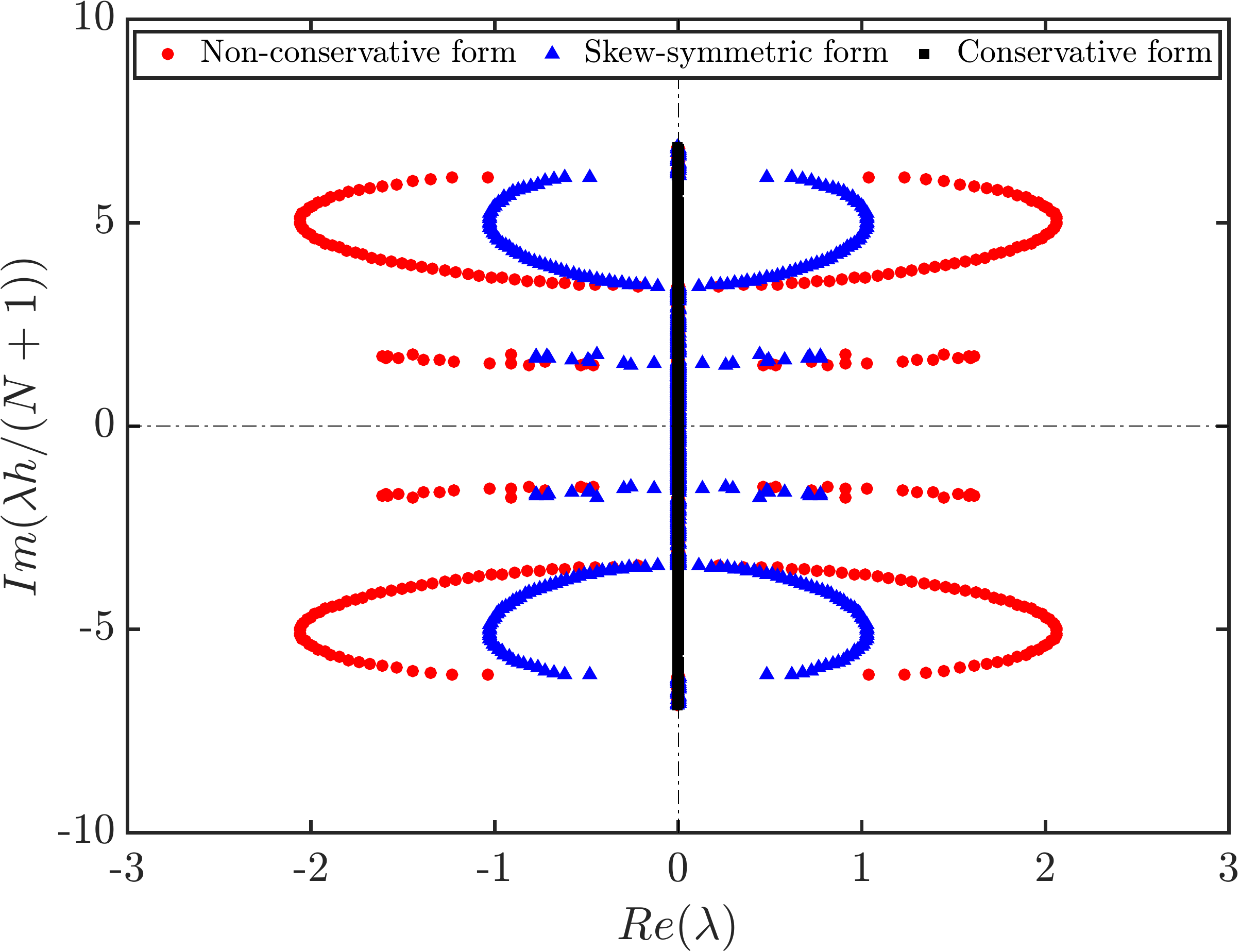}\label{fig:cons_spectra}}
	\subfigure[Non-conservative equation spectra]{\includegraphics[scale=0.23]{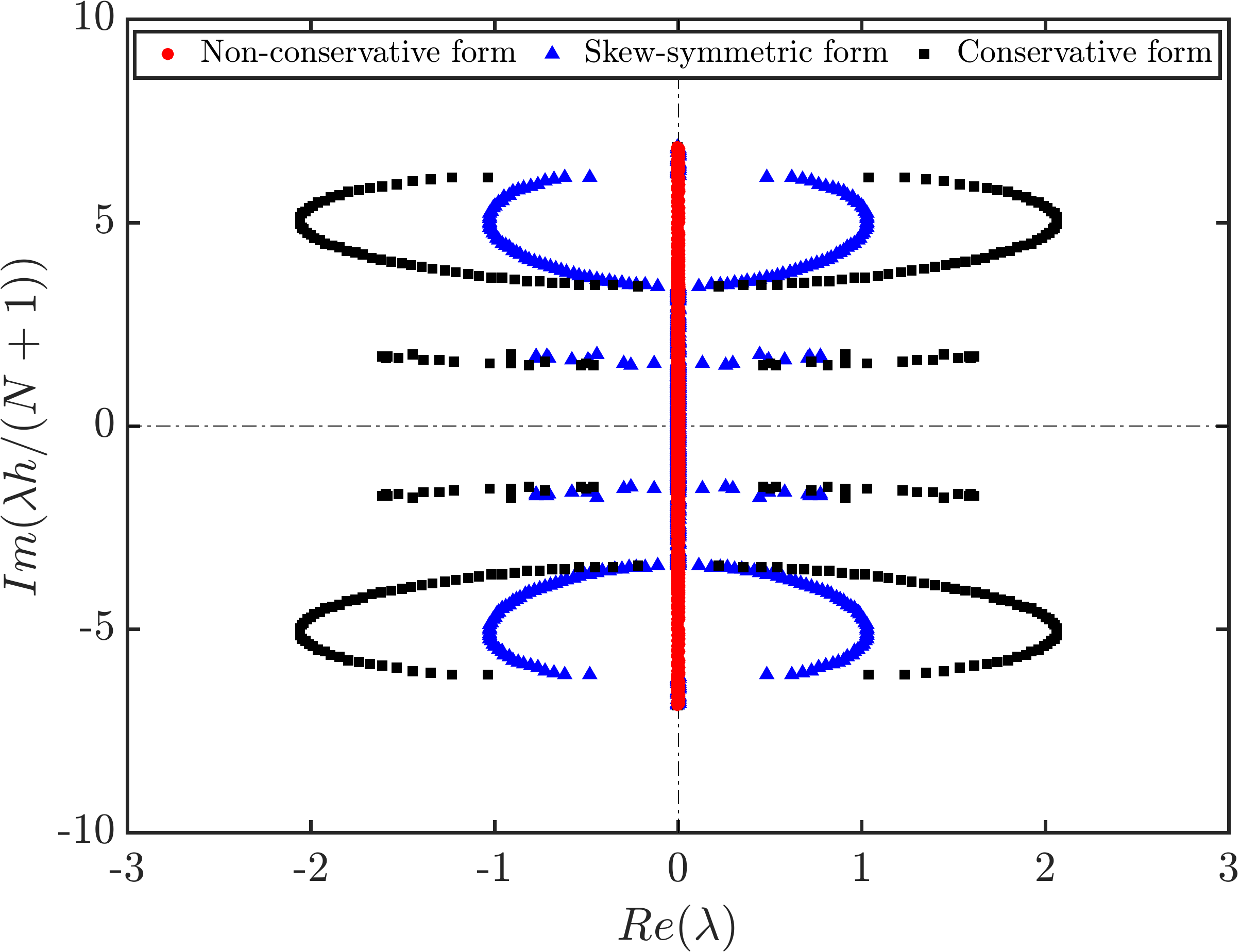}\label{fig:noncons_spectra}}	
	\caption{Eigenvalues of the non-constant speed \eqref{eq:a(x)} advection equation using Gauss-Lobatto points. The three relevant  versions (i.e. conservative ($\alpha=1$), skew-symmetric ($\alpha=1/2$), and non-conservative ($\alpha=0$) DG) of the split operator coefficient are represented. In this test case, $K=200$ elements have been used, while the polynomial order is $N=5$.}
\end{figure}

The eigenvalues for the conservative equation ($\theta=0$) are depicted in Figure~\ref{fig:cons_spectra}. The $x-$axes represents the real part of the spectra (where positive real parts lead to energy growth), whilst the $y-$axes represents their imaginary part. The imaginary part has been scaled with the element sizes and the polynomial order as in \cite{2015:Moura}. We see that both the non-conservative and skew-symmetric DG are unstable, whereas the conservative DG formulation is stable, consistent with the bound derived in \eqref{eq:conservative_L2_bound}. Note that the conservative DG spectra lies on the imaginary axis, and thus, none of the modes will exhibit energy growth or decay. On the other hand, when $\alpha\ne 1$, all modes are arranged by pairs, of which one of the pair shows exponential energy growth, and the other decay.\\

\begin{figure}
	\centering
	\subfigure[Conservative equation spectra]{\includegraphics[scale=0.23]{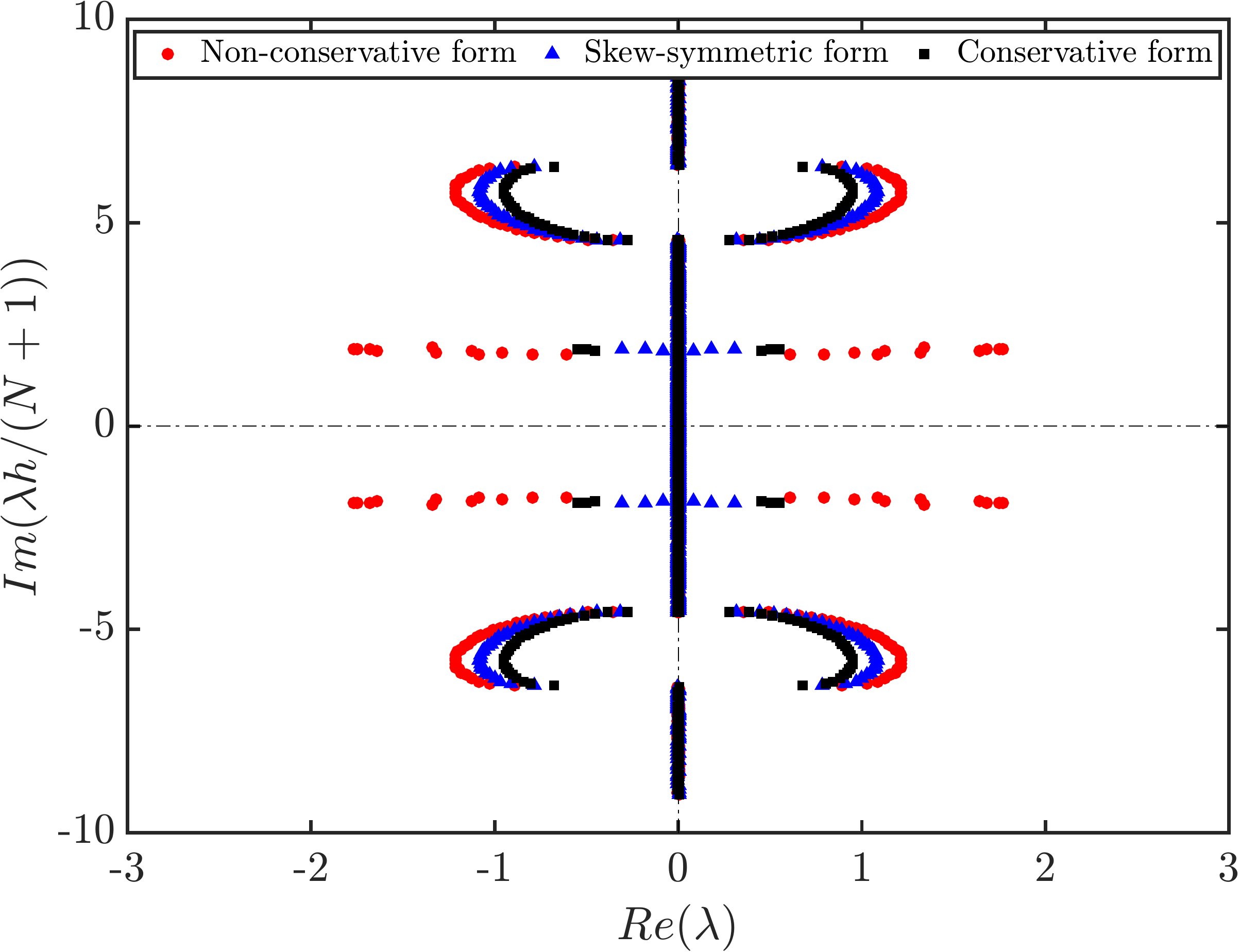}\label{fig:cons_spectra_LG}}
	\subfigure[Non-conservative equation spectra]{\includegraphics[scale=0.23]{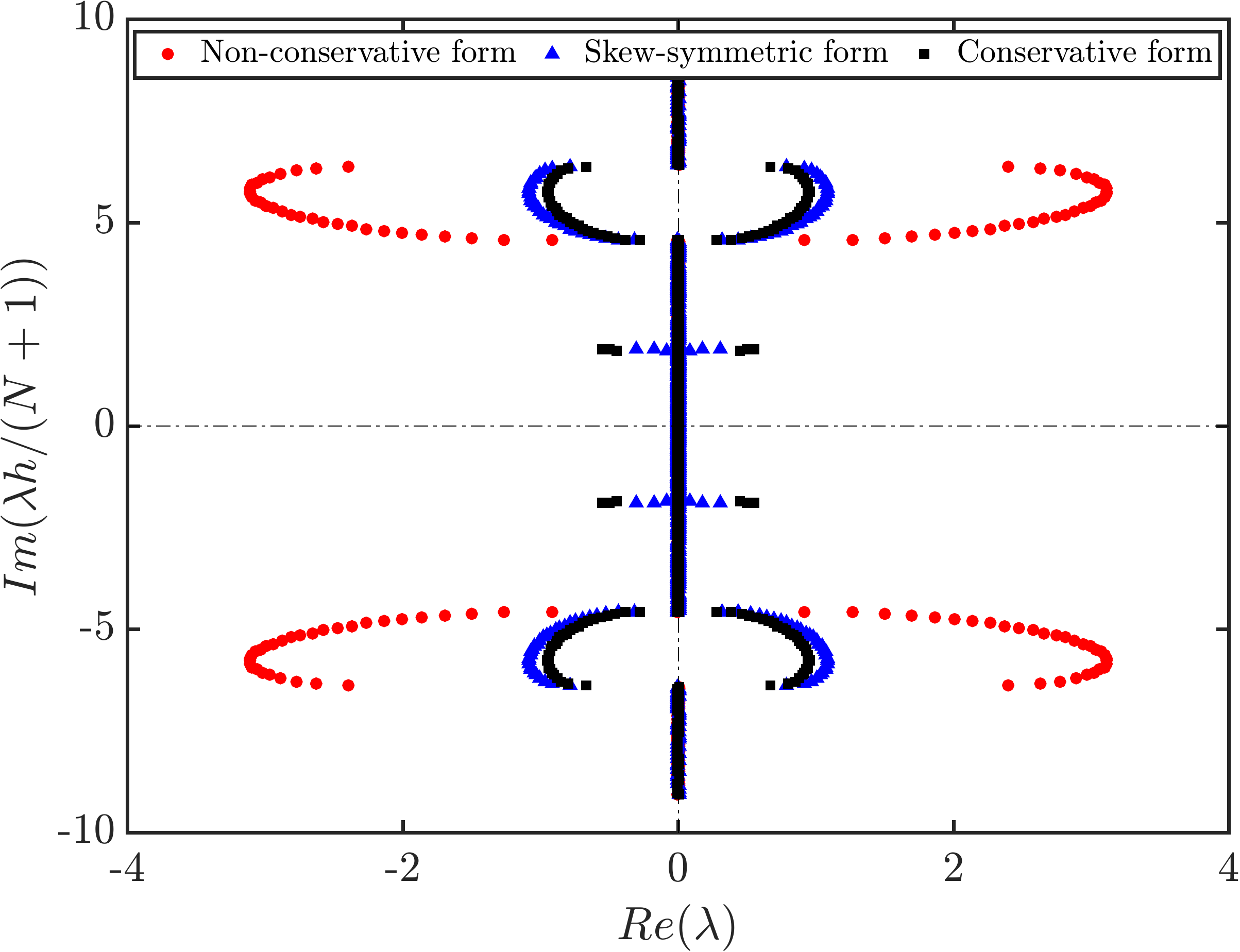}\label{fig:noncons_spectra_LG}}	
	\caption{Eigenvalues of the non-constant speed \eqref{eq:a(x)} advection equation using Gauss points. The three relevant  versions (i.e. conservative ($\alpha=1$), skew-symmetric ($\alpha=1/2$), and non-conservative ($\alpha=0$) DG) of the split operator coefficient are represented. In this test case, $K=200$ elements have been used, while the polynomial order is $N=5$.}
	\label{fig:spectra_LG}
\end{figure}
In Figure \ref{fig:noncons_spectra}, the non-conservative equation is analyzed. We observe that the results obtained are the opposite to the conservative equation ones, as shown by \eqref{eq:nonconservative_L2_bound}. In this case, the non-conservative DG discretisation is stable since its set of eigenvalues lie on the imaginary axes, which is consistent with the original PDE behaviour described in \eqref{eq:nonconservative_bound}.

To show the role played by the Gauss-Lobatto points, the spectra obtained using Gauss points is included in Figure \ref{fig:spectra_LG} for both the conservative and non-conservative problems. Recall that DG based on Gauss points does not satisfy the summation-by-parts property, and thus, the energy estimates that we have derived do not hold. An energy estimate for this problem was shown in \cite{2008:hest}. The computed eigenvalues show exponential energy growth for all the split operator approximations for both conservative (Figure \ref{fig:cons_spectra_LG}) and non-conservative (\ref{fig:noncons_spectra_LG}) problems. Hence, it does not seem to effectively remove aliasing errors by means of a split operator when using Gauss points. This energy growth (with Gauss points) must be then dissipated by other stabilisation techniques, such as adding artificial viscosity \cite{2006:Kirby}, or by means of interface dissipation with upwind Riemann solvers \cite{2009:Toro}. 

\begin{figure}
	\centering
	\includegraphics[scale = 0.3]{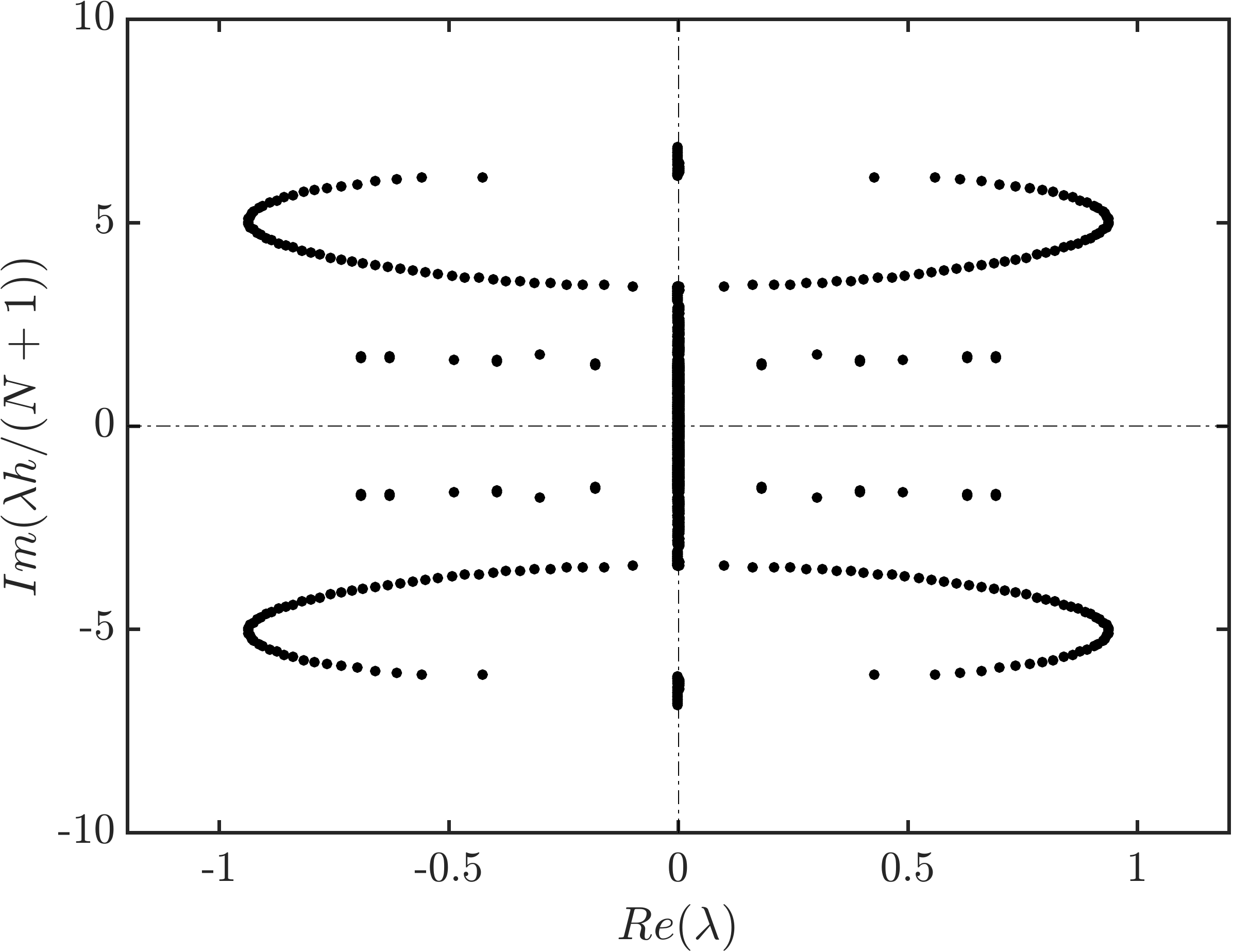}
	\caption{Eigenvalues of the conservative equation, solved with Gauss-Lobatto points and over-integration of the weak formulation integrals. The polynomial order is $N=5$, the quadrature order is $M=10$, and $K=200$ elements were used. Since the standard quadrature scheme yields a stable scheme, we can conclude that over-integration may lead to aliasing-driven instabilities. Its dissipation equals to that of the skew-symmetric formulation.}
	\label{fig:spectra_overintegration}
\end{figure}

Lastly, we show the effect of over-integration in the DGSEM-GL variant. We solve the conservative equation with the conservative DG (recall that when using over-integration, all split formulations are identical) since it was proven to be stable in its reduced (standard) quadrature version. However, when using over-integration, the scheme is not provably stable, according to \eqref{eq:bound_overintegration}. Figure \ref{fig:spectra_overintegration} depicts the numerical eigenvalues showing the aliasing driven instabilities arising as a result of the over-integration. The eigenvalues show the same growth/decay rates as the skew-symmetric form, since both follow the traditional bound, \eqref{eq:traditional_bound}. Therefore, over-integration, when considering the discontinuous Galerkin method with Gauss-Lobatto points must be used with care, since this technique does not necessarily stabilise the solution.

\section{Conclusions}\label{sec:conclusions}

In this work, we have studied the numerical instabilities that arise in the numerical solution of the non-constant advection speed equation. In particular, we have analysed non-constant advection speeds with constant sign. By means of energy estimates, we have shown that these instabilities are attributed to aliasing errors incurred in the numerical evaluation of the weak formulation integrals. Selecting the appropriate split form coefficient in a DG method with SBP-SAT properties is fundamental to achieve discrete energy conservation, whereas over-integration techniques are not capable of removing aliasing errors. Precisely, the conservative DG form satisfies these requirements when solving the conservative advection equation, whilst the same occurs when solving the non-conservative equation with the non-conservative DG version.  Discrete energy conservation does not occur when using DG versions that do not satisfy the SBP property, e.g. using the Gauss points.

\bibliographystyle{spmpsci}      
\bibliography{mybibfile}

\end{document}